\documentclass[12pt]{article}

\usepackage{amssymb,amscd,array}

\newcommand{\bbR}{\mathbb{R}}

\newcommand{\bbT}{\mathbb{T}}

\def\RP1{\mathop{\bbR P^1\!}\nolimits}

\newcommand{\Diff}{\mathrm{Diff}}

\def\Conf{\mathop{\rm Conf}\nolimits}

\def\PSL2{\mathop{\mathrm{PSL}(2,\bbR)}\nolimits}
\def\rg{\mathop{\textrm{\rm g}}\nolimits}

\newcommand{\cqfd}{\hspace*{\fill}\rule{3mm}{3mm}}

\begin{document}



\def\d{\delta}
\def\g{\gamma}
\def\om{\omega}
\def\r{\rho}
\def\a{\alpha}
\def\b{\beta}
\def\s{\sigma}
\def\vfi{\varphi}
\def\l{\lambda}
\def\m{\mu}
\def\implies{\Rightarrow}

\oddsidemargin .1truein
\newtheorem{thm}{Theorem}
\newtheorem{lem}[thm]{Lemma}
\newtheorem{cor}[thm]{Corollary}
\newtheorem{pro}[thm]{Proposition}
\newtheorem{ex}[thm]{Example}
\newtheorem{rmk}[thm]{Remark}
\newtheorem{defi}[thm]{Definition}


\title{Lorentzian worldlines and Schwarzian derivative}

\author{C.~Duval\footnote{mailto:duval@cpt.univ-mrs.fr}\\
{\small Universit\'e de la M\'editerran\'ee and CPT-CNRS}
\and
V.~Ovsienko\footnote{mailto:ovsienko@cpt.univ-mrs.fr}\\
{\small CNRS, Centre de Physique Th\'eorique}\thanks{CPT-CNRS, Luminy Case 907,
F--13288 Marseille, Cedex 9, FRANCE.
}
}

\date{}

\maketitle

\thispagestyle{empty}

The aim of this note is to relate the classical Schwarzian derivative
and the geometry of Lorentz surfaces of constant curvature.

\medskip
\noindent
{\bf 1.} 
The starting point of our investigations lies in the following remark
(joint work with L.~Guieu).
Consider a curve $y=f(x)$ in the Lorentz plane with metric $\rg=dxdy$.
If $f'(x)>0$, then its Lorentz curvature can be computed~:
$\varrho(x)=f''(x)\,(f'(x))^{-3/2}$
and enjoys the quite remarkable property~:
\begin{equation}
\label{rhoPrime}
\sqrt{f'(x)}\,\varrho'(x)
=
S(f)(x)
\end{equation}
where
\begin{equation}
\label{Schwarz}
S(f)(x)
=
\frac{f'''(x)}{f'(x)}-\frac{3}{2}\left(\frac{f''(x)}{f'(x)}\right)^2
\end{equation}
stands for the Schwarzian derivative of $f$. (It is well-known that the Schwarzian
derivative actually defines a quadratic differential we will write
$S(f)=S(f)(x)\,dx^2$.)

\medskip
\noindent
{\bf 2.} 
It is now natural to look for all Lorentz metrics admitting this specific
property.
We first consider an orientation preserving diffeomorphism
$f:\RP1\to\RP1$ whose graph is, therefore, a time-like curve of $\RP1\times\RP1$
endowed with the Lorentz metric $\rg=\rg(x,y)dxdy$, where $\rg(x,y)$ is a positive
function. Denoting by $t$ the Lorentz arc-length (also called proper time), we have

\begin{thm}
\label{Main}
The necessary and sufficient condition for which the equation
\begin{equation}
\label{drho}
d\varrho\,dt=S(f)
\end{equation}
holds true for any orientation preserving diffeomorphism $f$ of $\RP1$ is
\begin{equation}
\label{g}
\rg
=
\frac{dxdy}{(axy+bx+cy+d)^2}
\end{equation}
where $a,b,c,d$ are arbitrary real constants.
\end{thm}

Note that equation (\ref{drho}) is the intrinsic form of (\ref{rhoPrime}).
The metric (\ref{g}) is actually defined on the complement $\Sigma$ of the graph of
the linear-fractional transformation $y=-(bx+d)/(ax+c)$ associated with the singular
set of the metric. Clearly, $\Sigma$ has the topology of a cylinder $\bbR\times\bbT$.

\medskip
\noindent
{\bf 3.} 
The scalar curvature of the metric (\ref{g}) is constant, $R=8(ad-bc)$.
It is well-known \cite{Wolf,Kulk} that any Lorentz metric of constant curvature can
be locally brought into the following forms~:
\begin{eqnarray}
\label{g0}
\rg=dxdy\qquad\;\,& \hbox{if\ } R=0,\\[6pt]
\label{gR}
\rg=\displaystyle\frac{8}{R}\frac{dxdy}{(x-y)^2}& \hbox{if\ } R\not=0.
\end{eqnarray}
In our case (\ref{g}), this equivalence is global and can be obtained by the action
of $\PSL2\times\PSL2$ on $\RP1\times\RP1$. 

\medskip
\noindent
{\bf 4.} 
\textit{Proof of Theorem \ref{Main}}.
Recall that the curvature of a curve on a (pseudo-)Riemannian surface $(\Sigma,\rg)$
is given by $\varrho=\omega(v,a)\rg(v,v)^{-3/2}$, where $v$ stands for the velocity
and $a=\nabla_vv$ for the acceleration vector ($\omega$ is the surface $2$-form
associated with $\rg$ and $\nabla$ the Levi-Civita connection). For a time-like
curve $\tau\mapsto(x(\tau),y(\tau))$ of $\RP1\times\RP1$ one has
$$
\varrho
=
\frac{x'y''-x''y'}{\rg^{1/2}(x'y')^{3/2}}
-
\frac{x'\partial_x\rg-y'\partial_y\rg}{\rg^{3/2}(x'y')^{1/2}}
$$
as a function of the parameter $\tau$.
One then easily finds
\begin{equation}
\label{calculus}
\begin{array}{rcl}
\sqrt{\rg(v,v)}\,\varrho'
&=&
\displaystyle
-
\frac{x'''}{x'}+\frac{3}{2}\left(\frac{x''}{x'}\right)^2
+
\frac{y'''}{y'}-\frac{3}{2}\left(\frac{y''}{y'}\right)^2\\[10pt]
&&
\displaystyle
-
x'^2\left[
\frac{\partial_x^2\rg}{\rg}-\frac{3}{2}\left(\frac{\partial_x\rg}{\rg}\right)^2
\right]
+
y'^2\left[
\frac{\partial_y^2\rg}{\rg}-\frac{3}{2}\left(\frac{\partial_y\rg}{\rg}\right)^2
\right].
\end{array}
\end{equation}
In the r.h.s. of equation (\ref{calculus}) we recognize the
difference $S(y)(\tau)-S(x)(\tau)$ of the Schwarzian derivatives.

Let us show that
the extra terms vanish simultaneously iff the metric is given by (\ref{g}).
Indeed, if $\rg(x,y)=\partial_x\varphi(x,y)=\partial_y\tilde\varphi(x,y)$ with
$S(\varphi)(x)=S(\tilde\varphi)(y)=0$, then
$\varphi(x,y)=(\a(y)x+\b(y))/(\g(y)x+\d(y))$ and
$\tilde\varphi(x,y)=(\tilde\a(x)y+\tilde\b(x))/(\tilde\g(x)y+\tilde\d(x))$
with the unimodularity condition~: $\a\d-\b\g=\tilde\a\tilde\d-\tilde\b\tilde\g=1$.
Since
$\partial_x\varphi(x,y)
=1/(\g(y)x+\d(y))^2=
\partial_y\tilde\varphi(x,y)=1/(\tilde\g(x)y+\tilde\d(x))^2$,
the functions $\g,\d,\tilde\g$ and $\tilde\d$ turn out to be affine, whence
(\ref{g}).

Putting now $y=f(x)$ and $\tau=x$, and using the
definition of the arc-length~: $\rg(v,v)=f'(x)(dx/dt)^2=1$,
we readily get (\ref{drho}) from (\ref{calculus}).
\cqfd

\medskip
\goodbreak
\noindent
{\bf 5.} 
Amazingly, our standpoint allows us to recover the definition \cite{Hub} of
the relative Schwarzian derivative of two mappings of $\RP1$. 
\begin{cor}
Given a curve
$\tau\mapsto(x(\tau),y(\tau))$ of $\RP1\times\RP1$ as in Theorem \ref{Main}, the
equation (\ref{drho}) takes the form~:
\begin{equation}
\label{hubbard}
d\varrho\,dt
=
S(x,y)
\end{equation}
where $S(x,y)$ denotes the relative Schwarzian derivative of $x$ and $y$.
\end{cor}

\medskip
\noindent
{\bf 6.} 
Recall that in the classical Liouville theory associated with Riemannian
surfaces of constant curvature, the Schwarzian derivative enters naturally the
transformation law of the K\"alher metric under conformal transformations (see,
e.g., \cite{DNF} p. 118).

An analogous phenomenon occurs in the (real) Lorentz framework where the
Schwarzian derivative (\ref{Schwarz}) of a conformal diffeomorphism
$f\in\Conf(\Sigma,\rg)\cong\Diff(\RP1)$ for the metric (\ref{gR}) is interpreted in
\cite{KS} as the obstruction for $f$ to be an isometry. The conformal classes of the
metrics (\ref{g0},\ref{gR}) are studied in \cite{DG} where they are shown to be
symplectomorphic to coadjoint orbits of the Virasoro group. In these papers the
Schwarzian derivative appears as the
$1$-cocycle on the group $\Conf(\Sigma,\rg)$ and encodes the behavior of the metric
near the conformal boundary. 

Let us, however, emphasize that our results have no direct relationship with the
previous ones.

\medskip
\noindent
{\bf 7.}
Recently, E.~Ghys proved that \textit{given a diffeomorphism $f$ of $\RP1$, the
Schwarzian derivative $S(f)$ has at least four distinct zeroes} \cite{Ghys} (see
also \cite{OT,Tab}). This theorem was discovered as an analogue of the classical
four-vertex theorem~:
\textit{any smooth closed convex plane curve has at least four distinct curvature
extrema}. We have thus proved that the Ghys theorem is precisely the
four-vertex theorem for time-like closed curves in $\Sigma\subset\RP1\times\RP1$
endowed with the constant curvature metric (\ref{g}).

\medskip
\noindent
{\bf Acknowledgments.}
It is a pleasure to thank L.~Guieu for his precious help during the first stage of
this work and S.~Lazzarini for enlightening remarks.

\newpage


\end{document}